\RequirePackage{fix-cm}
\documentclass[twocolumn]{svjour3} 
\smartqed

\usepackage[utf8]{inputenc}
\usepackage{color}
\usepackage{graphicx}
\usepackage{xcolor}
\usepackage{amssymb}
\usepackage{lineno}
\usepackage{url}
\usepackage{bm}
\usepackage{framed}
\usepackage[]{natbib}
\journalname{Brain Topography}

\begin{document}
\title{Randomized Multiresolution Scanning in Focal and Fast   E/MEG Sensing of Brain Activity with a Variable Depth}

\author{A.\ Rezaei \and  A.\ Koulouri \and S.\  Pursiainen }
\institute{Atena Rezaei \at
               Faculty of Information Technology and Communication Sciences, Tampere University, P.O.\ Box 692, 33101 Tampere, Finland\\
              \email{atena.rezaei@tuni.fi} 
               \and 
            Alexandra Koulouri \at 
            Faculty of Information Technology and Communication Sciences, Tampere University, P.O.\ Box 692, 33101 Tampere, Finland \\
           \and   Sampsa Pursiainen \at 
             Faculty of Information Technology and Communication Sciences, Tampere University, P.O.\ Box 692, 33101 Tampere, Finland\\ }

\maketitle

\begin{abstract}

    We focus on elec\-tro-/mag\-ne\-to\-en\-cephalo\-graphy imaging of the neural activity and, in particular, finding a robust estimate for the primary current distribution via the hierarchical Bayesian model (HBM). Our aim is to develop a reasonably fast {\em maximum a posteriori} (MAP) estimation technique which would be applicable for both superficial and deep areas without specific {\em a priori} knowledge of the number or location of the activity. To enable source distinguishability for any depth,  we introduce a randomized multiresolution scanning (RAMUS) approach in which the MAP estimate of the brain activity is varied during the reconstruction process.  RAMUS aims to provide a robust and accurate imaging outcome for the whole brain, while maintaining the computational cost on an appropriate level. The inverse gamma (IG) distribution is applied as the primary hyperprior in order to achieve an optimal performance for the deep part of the brain. In this proof-of-the-concept study, we consider  the detection of simultaneous thalamic and somatosensory activity via numerically simulated data modeling the 14-20 ms post-stimulus somatosensory evoked potential and field response to  electrical wrist stimulation. Both a spherical and realistic model are utilized to analyze the source  reconstruction discrepancies. In the numerically examined case, RAMUS was observed to enhance the visibility of deep components and also marginalizing the random effects of the discretization and optimization without a remarkable  computation cost. A robust and accurate  MAP estimate for the primary current density was obtained in both superficial and deep parts of the brain. 
    
\end{abstract}


\section{Introduction}

    This study concentrates on electro-/mag\-ne\-to\-ence\-pha\-lo\-graphy (E/MEG) imaging of the brain activity  \citep{he2018electrophysiological}. The present focus is on the hierarchical Bayesian model (HBM)   \citep{calvetti2009,lucka2012} which allows one to find a focal and robust reconstruction by exploring a posterior probability distribution following from a conditionally Gaussian prior model. Our aim is, in particular, to develop a fast  {\em maximum a posteriori} (MAP) estimation technique which would be applicable for both superficial and  deep areas without additional {\em a priori} knowledge of the brain activity, such as physiological depth weighting  \citep{calvetti2015hierarchical,calvetti2018brain,homa2013bayesian}. While high-density measurements \citep{seeber2019subcortical} and advanced signal processing strategies \citep{pizzo2019deep} have recently been shown to be essential in distinguishing deep activity, this study focuses on the importance to reduce the random effects of the numerical discretization and optimization errors on the reconstruction process. 
    
We introduce a randomized multiresolution scanning (RAMUS) method in which the MAP estimate of the brain activity is refined gradually in the reconstruction procedure. RAMUS aims at reducing the random effects of the numerical discretization on the final estimate. It processes the well and ill-conditioned parts of the source space separately which has been suggested for ill-posed problems, e.g., in \citep{pursiainen2008_3,liu1995sensitivity,piana1997projected}.  A multiresolution decomposition provides an approximative split between detectable and undetectable parts for different source depths, as the maximal source localization accuracy varies strongly with respect to the depth  \citep{tarkiainen20033d,cuffin2001realistically,cuffin2001spherical,grover2016fundamental,wang2009relationship} with only the low resolution fluctuations being visible in the deep part of the brain \citep{pascual1999low,pascual1999review}. At each resolution level,  a  MAP estimate is evaluated via the iterative alternating sequential (IAS) algorithm and the inverse gamma (IG)  hyperprior which has been found to be advantageous for detecting deep activity \citep{calvetti2009}. 
    
The previous results suggest that HBM can find a focal solution deep in the head via the Markov  chain Monte Carlo (MCMC) sampling techniques, especially, if the activity can be constrained into a region of interest (ROI)  \citep{calvetti2009,lucka2012}. However, processing large data sets involving temporal measurement sequences with an advanced MCMC  approach without {\em a priori} knowledge of a ROI might be computationally too expensive for the practical use. Therefore, finding a robust and fast approach to distinguish activity reliably is crucial regarding the practical applications. In this proof-of-the-concept study, we consider the detection of simultaneous somatosensory and thalamic activity with numerically simulated data. This setup models the detection of the somatosensory evoked potentials and fields (SEP/F) in response to the electrical stimulation of the median nerve, particularly, thalamic (deep) P14/N14 and somatosensory (superficial) P20/N20 component peaked at 14 and 20 ms post-stimulus, respectively   \citep{buchner1988serial,buchner1995somatotopy, buchner1994preoperative,buchner1994source,haueisen2007identifying,attal2013assessment,fuchs1998improving}. 

In the numerical experiments, both a spherical and realistic model has been used to analyze the source reconstruction  discrepancies with RAMUS. The results suggest that a randomized set of decompositions \citep{mallat1989theory,clark1995multiresolution} is essential to marginalize out the possible modeling errors due to projecting the source space into different resolution levels which, again, is necessary in order to achieve the depth-invariance of the final MAP estimate. 


\section{Methodology}
\label{methods}

\subsection{Observation model}

For the EEG source modelling, we employ the finite element method and the current preserving H(div) approach  \citep{pursiainen2012b,pursiainen2016,miinalainen2019} in which the primary current distribution of the neural activity is assumed to have a square-integrable divergence $\vec{J}^P \in H(\hbox{div}) = \{\vec{w} | \nabla \cdot \vec{w} \in L^2 (\Omega) \}$ in the source space denoted by $\mathcal{S}$. 
The observation model 
is 
\begin{equation} {\bf y} = {\bf L} {\bf x} + {\bf n}, 
\end{equation} where ${\bf y}\in\mathbb{R}^m$ is the measurement vector, ${\bf L}\in\mathbb{R}^{m\times 3K}$ is the lead field matrix, ${\bf x}\in\mathbb{R}^{3K}$ is the unknown primary current distribution with $K$ denoting the total number of the source positions, and ${\bf n} \in \mathbb{R}^m$ is the measurement noise vector which is modelled as Gaussian random variable with zero mean and covariance matrix of the form ${\sigma}^2 {\bf I} \in \mathbb{R}^{m\times m}$. In this numerical study,  the diagonal covariance is  used for simplicity as it allows fixing the noise level with a single parameter, i.e., the standard deviation ${\sigma}$.  We refer to  $\mathbb{R}^{3K}$ as the source space $\mathcal{S}$ for the inverse problem of finding ${\bf x}$ given the data ${\bf y}$. The number of sources is three times the number of their positions, as each position is assumed to have three sources oriented along the Cartesian coordinate axes.


\subsection{Hierarchical Bayesian model}
\label{sec:hbm}

In the HBM framework, the prior of ${\bf x}$ is not fixed but random. It is  determined by the realization of the so-called  hyperparameter $\bm \theta$. The hyperparameter follows an {\em a priori} assumed distribution, i.e., the hyperprior. Consequently, the prior is a joint density given by 
$
p ( {\bf x}, {\bf \bm \theta}) \propto  p({\bf \bm \theta}) \,  p ({\bf x} \mid {\bf \bm \theta})
$ 
of ${\bf x}$ and  ${\bf \bm \theta}$. The conditional part of the prior  $p ( {\bf x} \mid {\bf \bm \theta})$ corresponds to a zero mean Gaussian density with a diagonal covariance matrix predicted by the hyperprior $p ({\bf \bm \theta})$. The hyperparameter ${\bm \theta}$ is of the same dimension as ${\bf x}$ with each entry  defining  the variance of its  respective entry in ${\bf x}$. The density of the hyperprior is long-tailed, implying that ${\bf x}$ is likely to be a sparse vector with only few nonzeros, which is advantageous for finding a focal reconstruction of the brain activity. As a hyperprior, one can use, e.g.,  the gamma (G) or  inverse gamma $\hbox{IG}({\bf \bm \theta} \mid \beta, \theta_0)$ density \citep{calvetti2009}, whose shape and scale are controlled by the parameters $\beta$ and $\theta_0$, respectively. IG is a  conjugate prior for a Gaussian distribution with an unknown variance (here the conditional prior), meaning that the corresponding posterior (here the actual prior) is also Gaussian. Again, G is a conjugate prior with respect to the reciprocal of the variance \citep{ohagan2004}.

The {\em posterior} probability density of ${\bf x}$, following from the classical Bayes formula \citep{ohagan2004}, is of the form
\begin{equation}
p ( {\bf x}, {\bf \bm \theta} \mid {\bf y})  =   \frac{p({\bf x}, {\bf \bm \theta}) \, p({\bf y} \mid {\bf x})} {p({\bf y})} \propto  {p({\bf x}, {\bf \bm \theta}) \, p({\bf y} \mid {\bf x})}, 
\end{equation} 
i.e., it is proportional to the product between the  prior density $p({\bf x}, {\bf \bm \theta})$, and the likelihood function $p({\bf y} \mid {\bf x}) \propto \exp ( - (2  \sigma^2)^{-1} \| {\bf L} {\bf x} - {\bf y}\|^2 )$ given by the measurement noise model \citep{schmidt1999bayesian}. 


We consider finding the inverse estimate via the {\em iterative alternating sequential} (IAS) MAP estimation method (Appendix \ref{sec:ias})  using primarily  the IG density as the hyperprior. IG has been suggested for depth localization in \cite{calvetti2009}, where the IG and G based IAS MAP estimate have been  shown to correspond to the minimum support and minimum current estimate (MSE and MCE)  \citep{nagarajan2006controlled}, respectively,  while the  first step of the iteration concides with the classical minimum norm estimate (MNE) \citep{hamalainen1993}. A recent comparison between IAS and other brain activity reconstruction techniques can be found in  \citep{calvetti2018brain}.  

The numerical exploration of the posterior density $p({\bf x}, {\bf \bm \theta} \mid {\bf y})$  is subject to the numerical discretization, i.e., the numerical definition of the source space $\mathcal{S}$ for ${\bf x}$ and the resulting lead field matrix. We aim to  reduce the effect of the discretization via the following two strategies motivating the  introduction of the RAMUS approach:
\begin{enumerate}
\item 
The reduction of the source space is essential to improve the ability of a solver to recover focal sources both in deep and superficial locations. Furthermore, since a sparse source space results here in source reconstruction of low spatial resolution, a source space refinement during the reconstruction process of this study is crucial. 
\item A  randomized set of decompositions 
enables averaging out (marginalizing) the effect of the
discretization error. 
\end{enumerate}
The theoretical justification of 1.\ and 2.\ are given in the following Sections \ref{sec:coarse_to_fine} and \ref{sec:randomized_scan}, respectively.

\subsubsection{Coarse-to-fine optimization}
\label{sec:coarse_to_fine}

\begin{figure*}[h!]
     \centering
    \includegraphics[width=10.5cm]{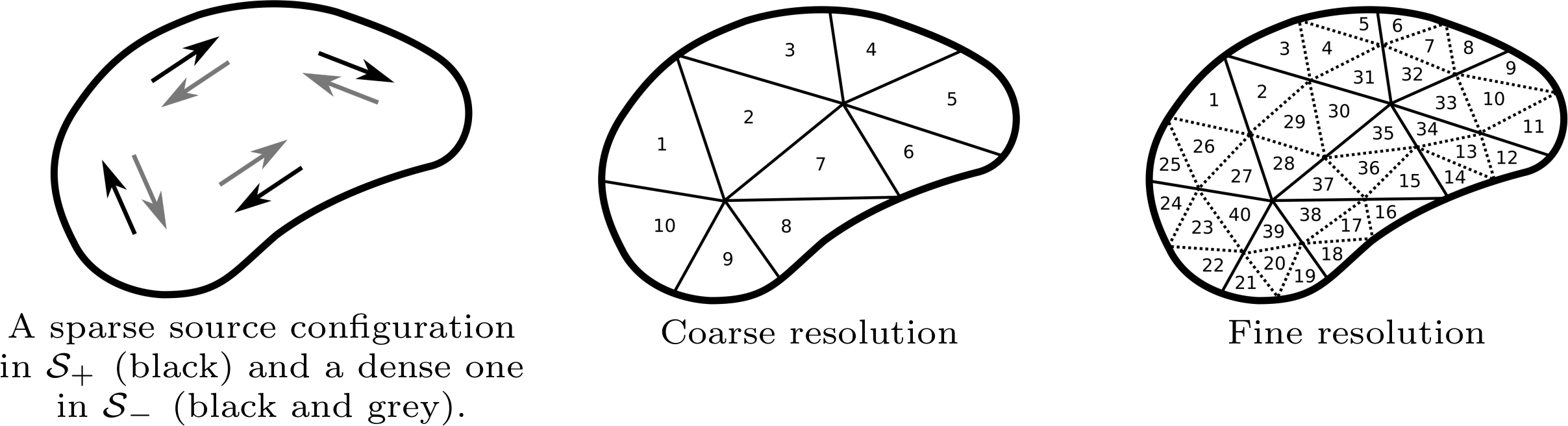}
    \caption{{\bf 1st from the left:} In E/MEG,  a coarse enough source configuration can be distinguished, i.e., it belongs to  $\mathcal{S}_\varepsilon^+$, while a dense one has modes that cancel each other and might be indetectable, i.e., in  $\mathcal{S}_\varepsilon^-$ (Figure \ref{fig:refinement}). {\bf 2nd and 3rd from the left:} An example of subdividing the grey matter compartment to subdomains in the case of a coarse (center) and fine (right) resolution. Here, the sparsity factor $s$, i.e., the ratio between number of subdomains for two consequtive resolution levels, would be  four. An example of mapping a subdomain from a coarse to fine resolution is given by $\{ 2 \} \to \{2,28,29,30\}$. }
    \label{fig:refinement}
\end{figure*}

\begin{figure}[h!]
 \centering
   \includegraphics[width=5.5cm]{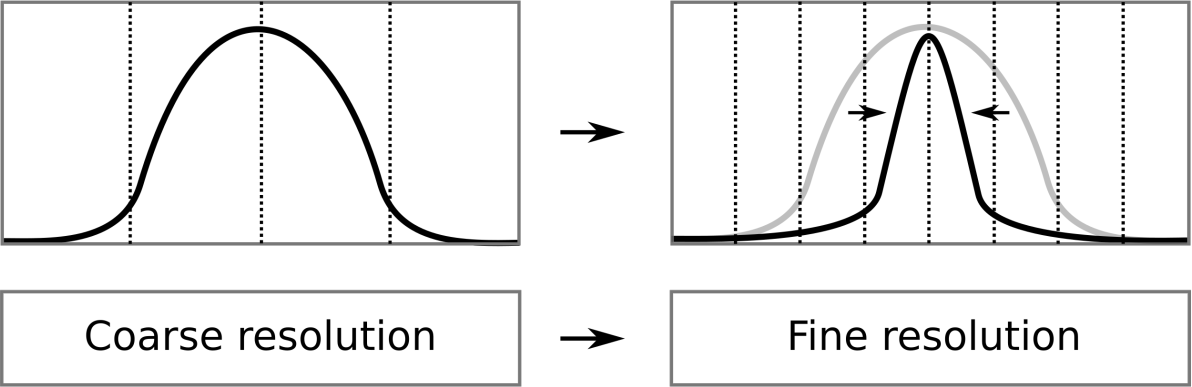}
     \label{fig:coarse_to_fine}

    \caption{
    Once an approximation for a non-zero source has been found at a coarse resolution level (left) the its support will shrink at the finer levels (right).\label{fig:coarse_to_fine}}
\end{figure}

The EEG source imaging problem is severely ill-posed \citep{Grech2008} and it is well-known that most of the solvers suffer from depth bias effects \citep{pascual1999review,Koulouri2017,Awan2018}. A way to reduce the ill-conditioning in the computations is by introducing coarser (sparse) source space, i.e.,  regularization by discretization \citep{hansen2010discrete,Kirsch2011}, or by approximating  the source distribution as a linear combination of
spatial basis functions (redundant
dictionaries) as proposed in
\cite{Haufe2008}. With a dimensionality reduction, the linear system to be solved is often over-determined and stable estimates can be obtained. However, this comes at a cost of poor resolution reconstructions due to large discretization errors. The idea of employing a multiresolution approach \citep{mallat1989theory},  where a progressive refinement in the source space is performed in order to obtain more accurate estimates, has been proposed for the E/MEG problem for example in
 \cite{gavit2001multiresolution,Malioutov2005}.

The source space $\mathcal{S}$ can be decomposed via the direct sum of  $\mathcal{S}_\varepsilon^+ = \{ {\bf 0}\} \cup \{ {\bf x} \, : \, \| {\bf L} {\bf x} \| \geq \varepsilon \}$ and $\mathcal{S}_\varepsilon^- = \{ {\bf x} \, : \, \| {\bf L} {\bf x} \| <  \varepsilon \}$, i.e.\ $\mathcal{S} = \mathcal{S}_\varepsilon^+ \oplus \mathcal{S}_\varepsilon^-$, where $\varepsilon$ is determined by the noise level.  $\mathcal{S}_\varepsilon^+$ and $\mathcal{S}_\varepsilon^-$ represent the sets of the detectable and undetectable source distributions, respectively. If possible, it is advantageous to decompose   $\mathcal{S}$ into $\mathcal{S}_\varepsilon^+$ and $\mathcal{S}_\varepsilon^-$ as, thereby, one can avoid source localization errors related to the indetectable distributions $\mathcal{S}_\varepsilon^-$  \citep{pursiainen2008_3,piana1997,liu1995sensitivity}. 
In E/MEG,  a coarse enough source configuration can be distinguished, i.e., it belongs to  $\mathcal{S}_\varepsilon^+$, while a dense one has modes that cancel each other and might be indetectable, i.e., in  $\mathcal{S}_\varepsilon^-$ (Figure \ref{fig:refinement}). The coarsity is specifically important considering deep activity for which the magnitude of the lead field is comparably low and, therefore, any deep source configuration is likely to belong to  $\mathcal{S}_\varepsilon^-$.  For a given lead field matrix  ${\bf L}$, the maximum possible number of detectable sources and, thereby, the maximal dimension of $\mathcal{S}_\varepsilon^+$ is  determined by the maximum number of nonzero singular values which coincides with the smaller dimension of ${\bf L}$, that is, the number of the data entries $m$. 

In the coarse-to-fine reconstruction strategy, the aim is to first limit the source space $\mathcal{S}$ to a subspace $\mathcal{S}_\varepsilon^+$ by restricting its resolution, to gradually increase its resolution, and  to eventually obtain an approximation for the whole space  $\mathcal{S}$. A nested set of restricted subspaces with different resolutions referred here to as a multiresolution decomposition is obtained recursively by selecting a uniformly random set of source positions from a given source space and associating those with the original set of positions through nearest interpolation. The coarsest resolution level is associated with the index $\ell = 1$. When moving from the $\ell$-th resolution level to the $(\ell+1)$-th one, the number of source positions is assumed to grow by a constant sparsity factor $s > 1$. An example for a dual resolution decomposition and a mapping of the subdomains between them can be found in Figure \ref{fig:refinement}.
 
In the IAS MAP estimation process, once the activity has been found at a coarse reconstruction level, the support of the candidate solution will shrink along with the increasing resolution (Figure \ref{fig:coarse_to_fine}). That is, the size of the details found is subject to the resolution level. Therefore, the final estimate is found as a combination of the estimates obtained for the different levels. In order to distinguish the weakly detectable activity, especially, the deep components, the number of the dimensions in the initial set should be of the same size with $m$, following from the maximal dimensionality of $\mathcal{S}_\varepsilon^+$.

\subsubsection{Randomized scanning}
\label{sec:randomized_scan}

Since a sparse source space is likely to induce a bias to the consequent estimates,  we propose to use a random set of (initial) sparse source spaces that aims to reduce the propagation of random discretization and optimization errors. The relationship between the global posterior optimizer   ${\bf x}^{\ast}$ and  ${\bf x}_k$ for the original source space $\mathcal{S}$ and its restriction    $\mathcal{S}_k$, respectively, can be modeled via the equation 
\begin{equation}
\label{aku_ankka}
    {\bf x}_k  = {\bf x}^{\ast} + {\bf d}_k + {\bf v}_k, 
\end{equation}
where ${\bf d}_k$ and ${\bf v}_k$ represent a discretization and optimization error, respectively. Of these, ${\bf v}_k$ depends of the quality of the MAP  optimization method and vanishes in the ideal case, while  ${\bf d}_k$  is fixed. If the degrees of freedom in $\mathcal{S}_1, \mathcal{S}_2, \ldots, \mathcal{S}_D$ have an independent and identical random distribution, the respective discretization errors ${\bf d}_1, {\bf d}_2, \ldots, {\bf d}_D$ can be modeled as independently and identically distributed random variables and, by the law of large numbers and the central limit theorem,  the discretization error term $\frac{1}{D} \sum_{k = 1}^D  {\bf d}_k$ of the mean \begin{equation} \frac{1}{D} \sum_{k = 1}^D  {\bf x}_k  = {\bf x}^\ast + \frac{1}{D} \sum_{k = 1}^D  {\bf d}_k + \frac{1}{D} \sum_{k = 1}^D  {\bf v}_k  \end{equation} is an asymptotically Gaussian variable  with expectation $\tilde{\bf d}$ and the rate of convergence $ \frac{1}{D} \sum_{k = 1}^D \to \tilde{\bf d}$  is of the order $O(D^{-1/2})$ with respect to the number of source spaces \citep{liu2001}. Consequently, the random effects of the discretization errors can be marginalized via estimating ${\bf x}^{\ast}$ in multiple randomly (independently and identically) generated source spaces.  The expectation $\tilde{\bf d}$ can be regarded as the remaining systematic discretization error which is specific to the set $\mathcal{S}_1, \mathcal{S}_2, \ldots, \mathcal{S}_D$, i.e., the resolution level, and is related, for example, to the relationship between the maximal achievable level of detail and the structure of the actual unknown brain activity.    
 
\begin{figure*}[h]
\centering
     \includegraphics[height=3.5cm]{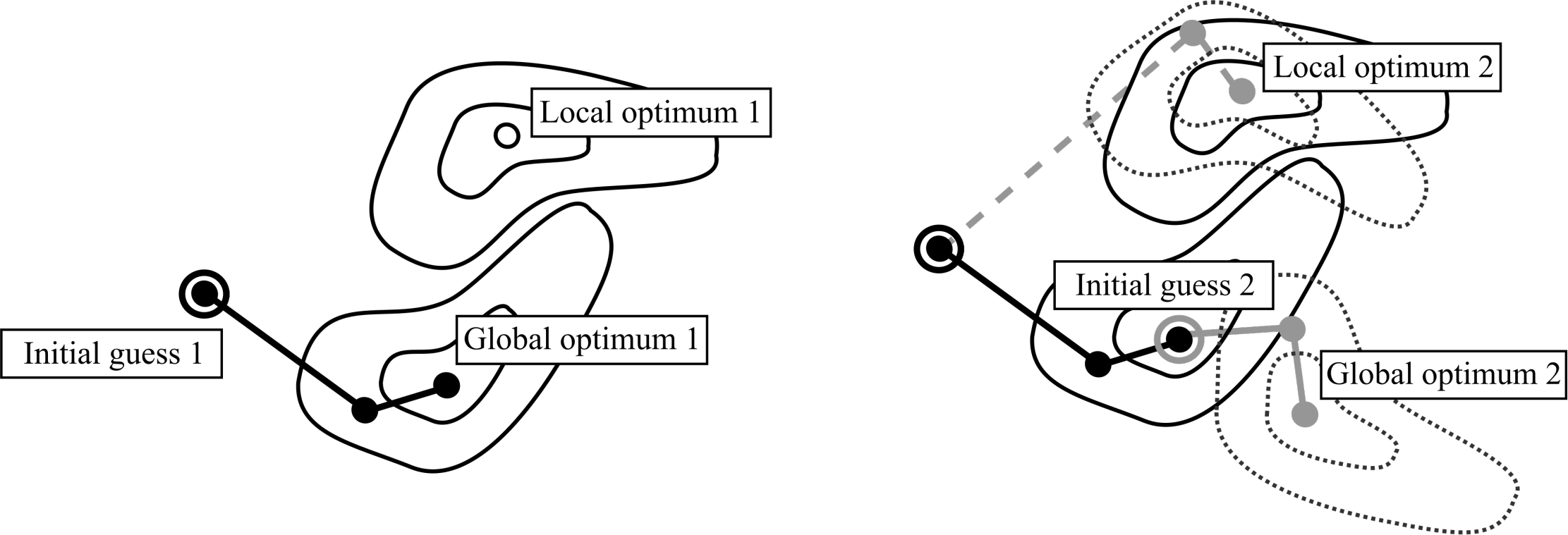} 
    \caption{An estimate for the global posterior optimizer  ${\bf x}_k$ obtained for the  source space $\mathcal{S}_k$ is found using the estimate for ${\bf x}_{k-1}$ as the initial guess (Section \ref{sec:randomized_scan}). We consider this approach necessary in order to maintain the estimates as close to the global optimum as possible.  Namely, using a fixed initial guess might mean that the global optimizer is not found for some of the source spaces. {\bf Left:} The global posterior optimizer is found for the posterior of space 1 (solid contours). {\bf Right:} For space 2 (dashed contours), it is found (solid grey path), if the final estimate obtained in the case one 1 is used as the initial guess for 2 (grey circled point), while a local optimizer is obtained (dashed grey path) with the original initial guess 1 (black circled point) resulting in an optimization error. The global optimum might, in practice, correspond to a situation in which both a superficial and deep source are detected, while the deep activity might be undetected at the local one.   }
    \label{fig:optimum_1_and_2}
\end{figure*}

Since the outcome of the optimization process for each given source space is {\em a priori} sensitive to the discretization errors, the estimate for ${\bf x}_k$ is found using the one for ${\bf x}_{k-1}$ as the initial guess. This approach is motivated by the present gradually progressing coarse-to-fine subdivision due to which the subsequent optimizers will be nearly similar. We consider it necessary in order to maintain each estimate in the vicinity of the  global optimum and, thereby, the norm of the optimization error ${\bf v}_k$ as small as possible.  Namely, using a fixed initial guess might mean that, instead  the global optimizer,  a local one is found for some of the source spaces as depicted in Figure \ref{fig:optimum_1_and_2}. The global optimum might correspond to a situation in which both a superficial and deep source are detected, while the deep activity might be undetected at a local one. 

Technically, updating the initial guess makes the estimate for ${\bf x}_k$ dependent on the previous one obtained for ${\bf x}_{k-1}$, i.e., the sequence of the estimates is  a time-homogeneous Markov chain. We regard the present approach as a surrogate transition rule \citep{liu2001} estimating the outcome of an ideal optimization method which would find the global optimum precisely with ${\bf v}_k = 0$, thereby, resulting in the identity
\begin{equation} \frac{1}{D} \sum_{k = 1}^D  {\bf x}_k  = {\bf x}^\ast + \frac{1}{D} \sum_{k = 1}^D  {\bf d}_k \to {\bf x}^\ast + \tilde{\bf d},  \label{average_approx} \end{equation} 
which will hold approximately, if the surrogate rule is accurate enough. 

\subsection{RAMUS}

\label{sec:RAMUS}

We propose the following algorithm for RAMUS to reduce the random discretization and optimization effects when finding a reconstruction for the unknown parameter ${\bf x}$ with the IAS MAP estimation method. 
\begin{enumerate}
\item Choose the desired number of the resolution levels $L$ and the sparsity factor (the ratio of source counts) $s$ between each level. The number of the sources at a given resolution level will be $K_\ell = K s^{(\ell - L) }$, where $\ell = 1, 2, \ldots, L$ is the index of the resolution level, the larger the value of the index $\ell$ the finer the resolution.
\item For each resolution level $\ell = 1, 2, \ldots, L$, create a random uniformly distributed set of center points $\vec{p}_1, \vec{p}_2, \ldots, \vec{p}_{K_\ell}$. Find source point subsets $B_1$, $B_2$, $\ldots$, $B_{K_\ell}$ applying the nearest interpolation scheme with respect to the center points. That is, each subset $B_j$ consists of those source positions of the total source space $\mathcal{S}$, whose nearest neighbor with respect to $\vec{p}_1, \vec{p}_2, \ldots, \vec{p}_{K_\ell}$ is $\vec{p}_j$. The  average number of source positions associated with $B_j$ is approximately given by the sparsity factor $s$. The resolution of this subdivision grows along the number of the center points. The unknown parameter is assumed to be constant in each subset, and the actual source count is assumed to stay unchanged regardless of the resolution.  
\item Repeat the first two steps to generate a desired number $D$ of independent multiresolution decompositions $\mathfrak{S}_1$, $\mathfrak{S}_2$, $\ldots$, $\mathfrak{S}_D$. 
\item  Start the reconstruction process with the decomposition $\mathfrak{S}_1$ and a suitably chosen  initial guess ${\bf x}^{(0)}$.
\item  For decomposition $\mathfrak{S}_k$,  find a reconstruction ${\bf x}^{(\ell)}$ with the IAS MAP technique with the initial guess ${\bf x}^{(\ell-1)}$ for the resolution levels $\ell = 1, 2, \ldots, L$.  
\item After going through all the decompositions, obtain the final estimate for the decomposition (basis) $k$ as the normalized mean 
\begin{equation}  \overline{\bf x}^{(k)} =   \sum_{\ell = 1}^L {\bf x}^{(\ell)} \,  /  \,  \sum_{\ell = 1}^L  s^{(L - \ell)}  , \label{average_1} \end{equation} 
where the denominator follows from the need to balance out the effect of the multiplied source count following from the interpolation of a coarse level estimate to a denser resolution level. 
\item If $k < D$, move to the next decomposition, i.e., update $k \to k + 1$, and repeat the previous step with the initial guess $\overline{\bf x}^{(k-1)}$ for the resolution level $\ell = 1$.
\item Obtain the final reconstruction as the mean: 
\begin{equation}
\overline{\overline{\bf x}}^{(k)} = \frac{1}{D} \sum_{k=1}^D \overline{\bf x}^{(k)}.  
\label{average_2}
    \end{equation}
\end{enumerate} 
Technically, this process is equivalent to first evaluating the mean (\ref{average_2}) for each resolution level and then the normalized mean (\ref{average_1}) over the different resolutions, showing that an approximation of the form (\ref{aku_ankka}) is, in fact, obtained for each set of independent and identically generated source spaces. Since the final reconstruction is obtained as a mean over all the reconstruction levels, also the potential systematic discretization errors will be averaged with an equal weighting. This approach is used, as different resolution levels localize different details (Section \ref{sec:coarse_to_fine}). Consequently, the details  found for the most levels are likely to gain the highest intensity in the final reconstruction. A schematic illustration of the resulting data flow has been included in Figure \ref{fig:update_scheme}.

\begin{figure}[h]
\centering
    \includegraphics[width=7.5cm]{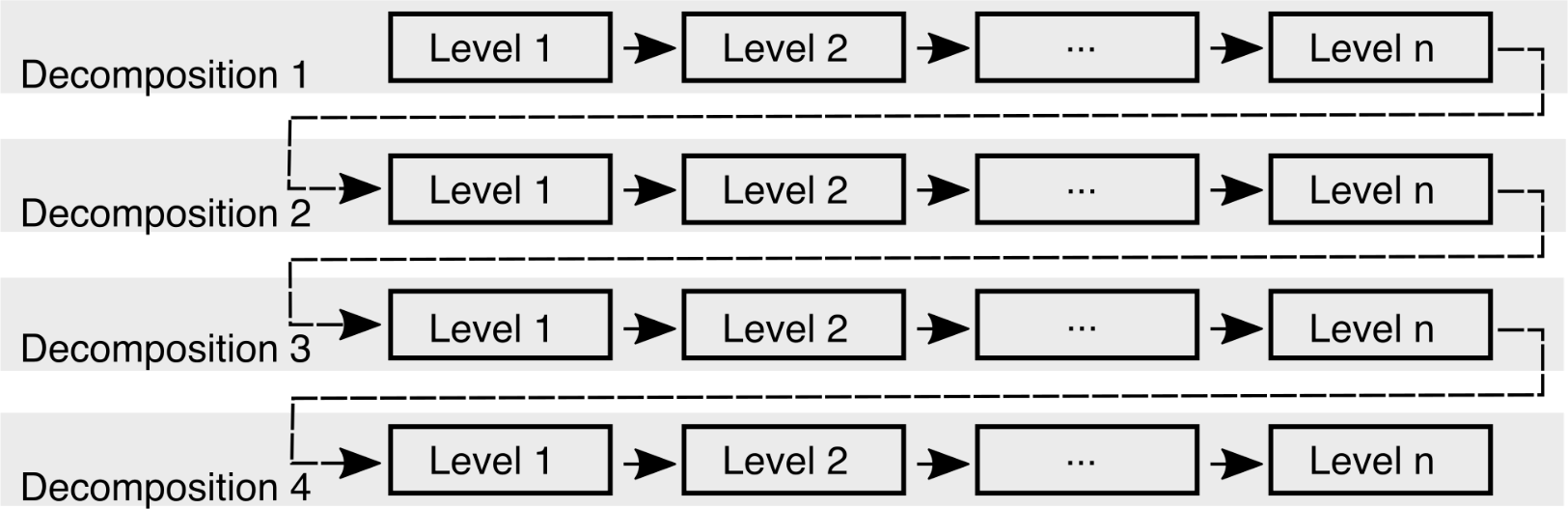}
    \label{fig:update_scheme}
    \caption{A schematic visualization of the data flow during the reconstruction process for the multiresolution decompositions $1, 2, \ldots, D$ each one with resolution levels $1, 2, \ldots, L$. The final estimate (\ref{average_1}) obtained for the decomposition $k$ is used as the level-one initial guess for the decomposition $k+1$. This sequential strategy for selecting the initial guess aims to minimize the effect of the optimization errors as suggested in Figure \ref{fig:optimum_1_and_2}. Note that with a good enough initial guess the global optimum is always found, meaning that the differences between the optimization results can be associated with the discretization errors which are modeled here as independently and identically distributed random variables. }
    \label{fig:update_scheme}
\end{figure}

\subsection{Numerical implementation with Zeffiro Interface}

The forward and inverse solvers applied in this study were implemented in the Matlab (The MathWorks Inc.) as a part of the {\em Zeffiro Interface} (ZI) code package which is openly available in GitHub\footnote{\url{https://github.com/sampsapursiainen/zeffiro_interface}}. ZI is a tool enabling finite element (FE) based forward and inverse computations in electromagnetic brain applications. The forward approach of ZI together with the basic version of the IAS source reconstruction approach have been validated numerically in \citep{miinalainen2019,pursiainen2012raviart}. ZI  generates a uniform tetrahedral finite element (FE) mesh. Each source distribution is obtained by picking the first $K$ entries of the randomly (uniformly) permuted set of the tetrahedron centers for the brain compartment. Due to the uniform mesh structure, this strategy leads to an evenly distributed set of source points. 
$\bf x$. 

ZI allows performing the source reconstruction routines using either a CPU or a GPU (graphics processing unit) type processor. Today, effective GPUs are available in power PCs an workstations but most laptops are still limited to CPU processing. Therefore, to compare the performance difference between GPU and CPU platforms, the computing time for forming a random set of multiresolution decompositions and inverting a given measurement data vector were evaluated for NVIDIA Quadro P6000 workstation GPU and Intel i7 5650U laptop CPU. 

\subsection{Numerical experiments}

In the numerical experiments, we used the realistic  population head model\footnote{\url{https://itis.swiss/virtual-population/regional-human-models/phm-repository/}} (PHM) \citep{lee2016investigational}, consisting of five layers (white matter, grey matter, cerebrospinal fluid (CSF), skull, and skin) and the three-layer Ary model in which concentric 87, 92 and 100 mm spheres present grey matter, skull and scalp layer. The  cerebellum and vetricle layers included in the PHM were modeled as part of the grey matter and CSF, respectively. The conductivity of each layer can be found in Table \ref{tab:conductivities}. The PHM and Ary model were discretized  with a uniform point lattice with the resolution 0.85 and  1 mm, leading to 24M and 30M tetrahedral elements and 4M and 5M nodes, respectively. In both cases, a single lead field matrix was generated for 10000 randomly distributed synthetic source positions. The lead field matrix entries were  evaluated in SI units, i.e.,  Ohm/m and 1/m\textsuperscript{2} for EEG and MEG, respectively. Each point contained  three sources oriented along x-, y- and z-directions. Since the grey matter compartment of PHM does not include the thalamus, the source space was extended to cover both the white and grey matter compartment. Note that the lead field matrix and the corresponding source space have to be generated only once after which the space can be decomposed in multiple ways, e.g., different resolutions, as is the case in the proposed RAMUS process.

\mbox{} 

\begin{table}[!h]
    \centering
        \caption{The conductivities (S/m) of the different compartments for PHM and Ary model \citep{ary1981}. Justification of the values used for the realistic model can be found in \citep{dannhauer2010}.}
    \label{tab:conductivities}
    \begin{footnotesize}
    \begin{tabular}{l r r r r r r}
 Model & WM & GM & CSF & Skull & Skin   \\ 
 \hline 
Ary  & &  0.33 & & 0.0042 & 0.33   \\ 
Deep & 0.14  & 0.33  &  1.79 & 0.0064 & 0.43 \\
    \end{tabular}
    \end{footnotesize}
\end{table}

\subsubsection{Simulated measurements}

\begin{figure}[h!]

 \centering
     
    \includegraphics[width=6cm]{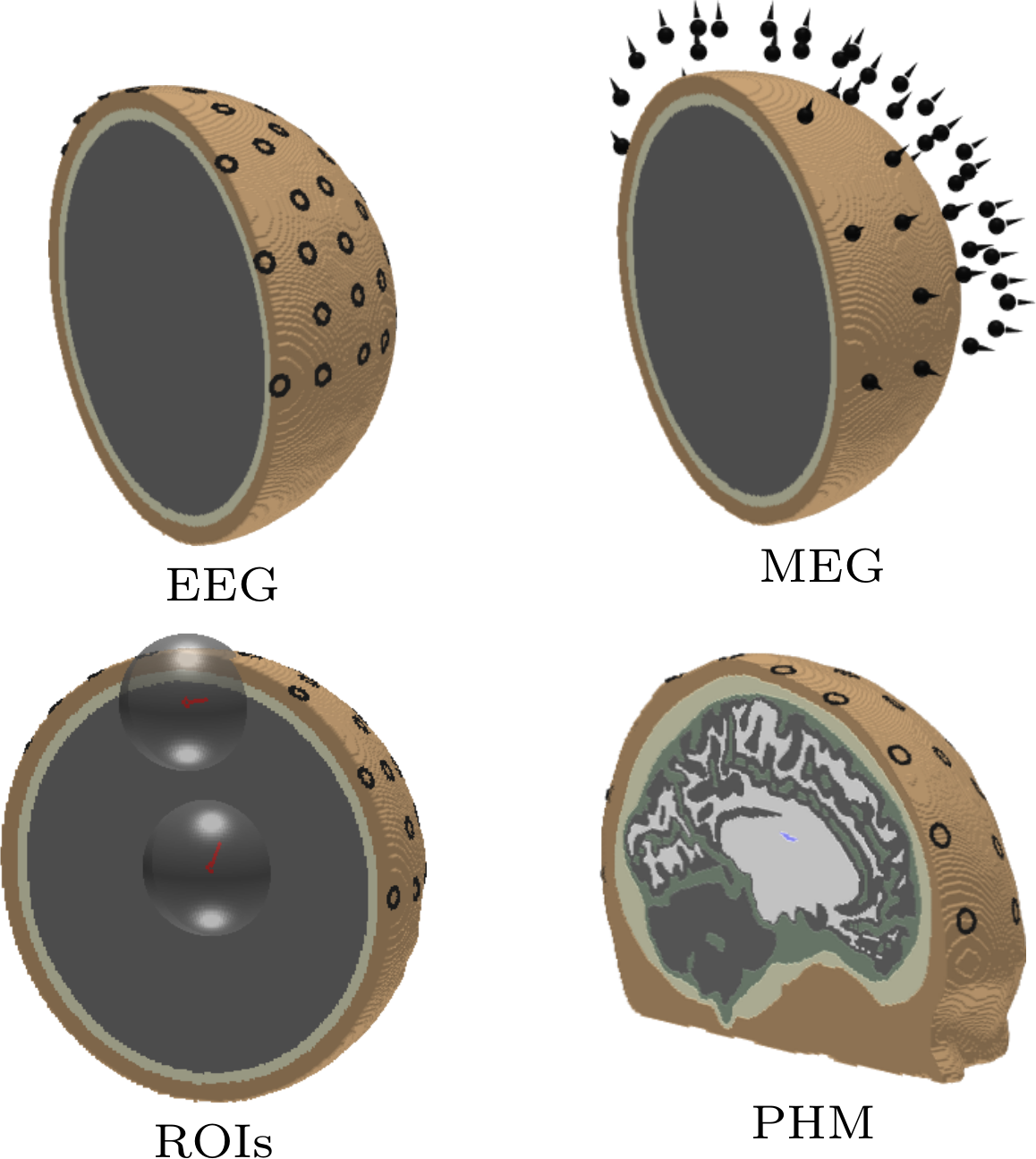}
    \caption{The volumetric FE mesh for the realistic five-layer PHM and three-layer Ary model. {\bf Top row:} The left image shows the domain with 102 EEG ring electrodes on it and the right one visualizes the positioning of the 102 radial magnetometers which has been obtained by scaling the electrode locations by the factor of 1.2. {\bf Bottom row:} The source locations with the ROIs for the Ary model (left) and a cut-view of the PHM (right) with 1cm diameter ring electrodes which were modeled using the complete electrode model \citep{pursiainen2012}. }
    \label{fig:domains}
\end{figure}

For the Ary model, a total of 102 sensor points were distributed over the upper hemisphere. Using those, both electrode and radial magnetometer measurements of the electromagnetic field were simulated as shown in Figure \ref{fig:domains}. The magnetometer locations were obtained by scaling the radial component of the source locations by a factor of 1.2. The electrodes were modeled using the complete electrode model \citep{pursiainen2012}. The inner and outer radius of the ring were 5 and 10 mm, respectively. The average contact resistance of each electrode was assumed to be 1 kOhm. In the case of PHM, an EEG cap with 72 ring electrodes (10 mm outer and 5 mm inner diameter, 1 kOhm resistance) was attached on the head model. 

Two current dipoles were placed in shallow and deep parts of the grey matter. The source locations can be found in Table  \ref{tab:source_specs}. Physiologically these could be interpreted as the somatosensory (superficial) P20/N20 and thalamic (deep) P14/N14 component, i.e., the 20 and 14 ms post-stimulus peaks. Activity for both locations occurs at the same time in the SEP/F response to the median nerve stimulus  \citep{buchner1988serial,buchner1994preoperative,buchner1995somatotopy}. When active simultaneously, the deep source was assumed to be slightly stronger in magnitude compared to the superficial one to enable the visibility of the deep part. This situation occurs momentarily in the median nerve stimulation, since the thalamic  source obtains its maximum before the somatosensory activity increases in magnitude.

 As the measurement error term, we used zero mean Gaussian white noise with standard deviation of 3\% respect to maximal signal amplitude. To investigate the noise-robustness of the source reconstruction, 5\% noise was used in a single test. For the generality of the results, the maximum data entry of each dataset was normalized to one. The  accuracy of the source recovery was analyzed in two 60 mm diameter spherical ROIs centered at the source locations (Figure  \ref{fig:domains}). 
\mbox{} 

\begin{table}[!h]
    \centering
        \caption{The source locations and orientations utilized in the numerical experiments.}
        \begin{footnotesize}
    \label{tab:source_specs}
    \begin{tabular}{l r r r r r r}
 Type & Corresp.\ &  x  & y  & z  & Angle \\ 
  & &(mm) &(mm) &(mm) & (deg) \\ 
 \hline 
Superf.\ & P20/N20 & -5 & 0 &  77 & 11  \\ 
Deep & P14/N14  & 7   &  0 & 5 & 68  \\
    \end{tabular}
    \end{footnotesize}
\end{table}

\subsubsection{IAS MAP iteration}

The previous experience shows that, in order to distinguish deep activity \citep{calvetti2009}, the hyperparameter values for the hyperprior have to be set as small as possible without risking the numerical stability of the reconstruction process. In the present study, the scale parameter  $\theta_0$ was chosen to be 1E-10 and the shape parameter $\beta$ was given the smallest possible value  1.5. These values were found to work generally well and they are supported also by the earlier studies \citep{calvetti2009}. Ten iteration steps were performed to obtain a MAP estimate for a single resolution level. A single step was utilized in a single test.

\subsubsection{Validation tests}

We analyzed the performance of the RAMUS reconstruction approach both visually and numerically in the  tests (A)--(I) using the Ary model. The spherical domain was used in order to optimize the clarity of the results. In addition to these reconstructions, one test (J) was performed using PHM, i.e., the realistic model. The specifications for (A)--(J) can be found in Table \ref{tab:A-H}.  

The accuracy obtained in the cases (A)--(I) was analyzed by comparing the average position (center of mass), orientation and magnitude of the reconstructed distribution within the ROI to that of the actual dipole source. These average estimates were obtained with respect to the final reconstructed distribution of 10000 sources in each case (A)--(I) and for both single and multiple resolution reconstructions. In addition, the relative magnitude (between 0 and 1) of the distribution was calculated for each ROI. The source was classified  as detected, if the relative maximum exceeded the value 0.1, and otherwise undetected.  This threshold criterion was chosen as it represents roughly the limit of a visually detectable source.  In (A)--(I), we varied the number of multiresolution decompositions, sparsity factor, hyperprior, the source magnitudes, and the measurement modality (EEG or E/MEG). When combining the lead field matrices for E/MEG, the MEG lead field matrix and data was scaled so that the Frobenius norm, i.e., the 2-norm of all the entries, was equal to that of the EEG lead field matrix. 

The case (I), was  studied using three alternative approaches in addition to the basic multiresolution scheme. In the first one of these, the noise level was increased to 5\%. The second one involved only the coarse resolution level with otherwise unchanged parameters. In the last one, only single IAS MAP iteration was performed on each reconstruction level, meaning that the estimate obtained coincided with MNE \citep{calvetti2009}.

\mbox{} 

\begin{table*}[!h]
        \caption{The specifications of the reconstructions computed in the numerical experiments.}
        \begin{footnotesize}
    \label{tab:A-H}
      \begin{center}
    \begin{tabular}{l l l r r r r r }
       & & & &  &  &  \multicolumn{2}{c}{Amplitude} \\ ID  & Geom. & Data & $s$$\mbox{}^a$ & Dec.$\mbox{}^b$ & Hyperprior   &Deep & Superf. \\ \hline 
        (A) & Ball & {\footnotesize  EEG} & 8 &  100 & IG  &  10 & 5 \\
        (B) & Ball &{\footnotesize  EEG} & 8  &  100 & IG  &  10 & 0 \\
        (C) & Ball &{\footnotesize  EEG} & 8 & 100 & IG &  0 & 5 \\
        (D) & Ball &{\footnotesize  EEG} & 8  & 100  & IG &  10 & 7 \\
        (E) & Ball &{\footnotesize  EEG} & 5 & 100  & IG & 10 & 7 \\
        (F) & Ball &{\footnotesize  EEG} & 8 & 20  & IG &  10 & 7 \\
        (G) & Ball &{\footnotesize  EEG} & 8 & 100  & G &  10 & 5 \\
        (H) & Ball &{\footnotesize  E/MEG} & 8 &  100  & IG &  10 & 5 \\
        (I) & Ball &{\footnotesize  E/MEG} & 8 & 100 & IG & 10 & 7 \\
        (J) & PHM & {\footnotesize  EEG} & 8 &  100  & IG &  10 & 7 \\
    \end{tabular}
    \end{center} 
    $\mbox{}^a$ Scaling factor. \\
    $\mbox{}^b$ Number of decompositions. \\
    \end{footnotesize}
\end{table*}

A total of 100--400 source positions at the coarse level, i.e., a number roughly comparable to that of the data entries (Section \ref{sec:coarse_to_fine}), was found to work appropriately in the detection of the deep activity. When the sparsity factor between $s=8$ and $s=5$, the source position count was within this interval at the coarsest level of a three-level multiresolution decomposition for the initial set of 10000 source positions. At the coarsest level, each source position was associated to about $s^2$  ($s^{L-1}$ with $L = 3$), i.e., between 64 and 25 finest-level source positions, respectively. The number $D$ of multiresolution decompositions was chosen to be comparable to this number, slightly below or above that, in order to guarantee sufficient averaging over all the possible random basis choices. 

\section{Results}

\begin{figure}[h!]
     \centering
   
    \includegraphics[width=6.5cm]{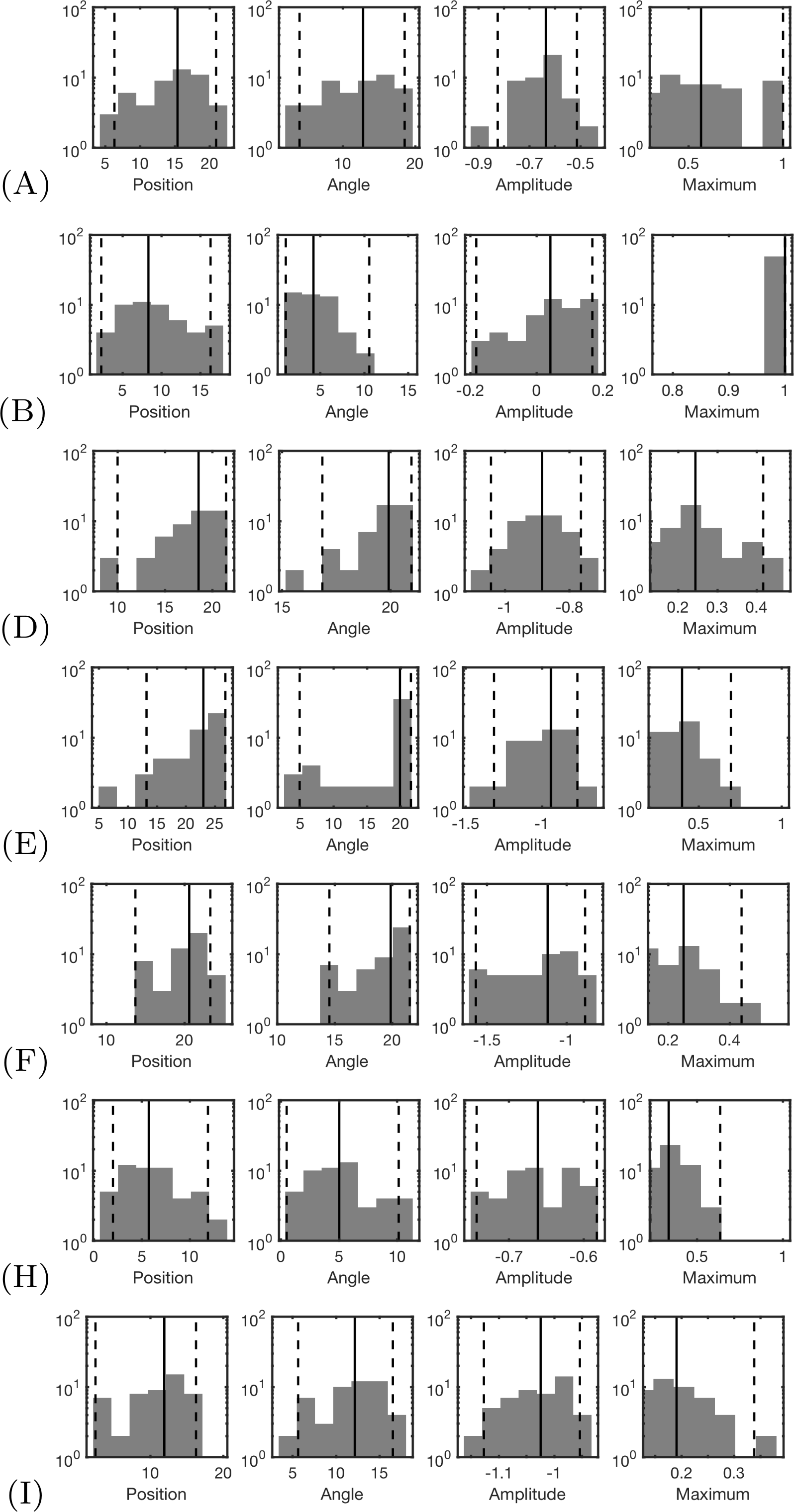}
    \caption{The results of the deep source localization for the numerical experiments (A)--(I) conducted in the spherical domain (Table \ref{tab:A-H}). The distributions of the position (mm), angle (deg) and relative logarithmic (log10) amplitude  difference to the exact dipole source, computed in the ROI, have been analyzed as histograms. The sample size is 50. Each reconstruction in the sample has been obtained by reconstructing the activity in the whole brain for an independent random realization of the noise vector and associating the total integrated activity in each ROI to the corresponding (deep or superficial) dipole source. Additionally, the histogram of the  relative maximum in the ROI is given. The solid vertical line shows the median for each distribution, and the dashed lines  mark the 90\% confidence interval. In general, the results show that the IG hyperprior is necessary for detecting the deep  source. The accuracy and reliability of the results increase along with the number of multiresolution decompositions. Furthermore, using E/MEG instead of EEG increased the accuracy of the deep source localization, while EEG was advantageous with respect to the amplitude of the deep source. The results are not  visualized for the cases in which the localization criterion (relative maximum $> 0.05$) was satisfied by less than 5\% of the reconstructions.}
    \label{fig:histograms_deep}
\end{figure}

\begin{figure}[h!]
     \centering
    \includegraphics[width=6.5cm]{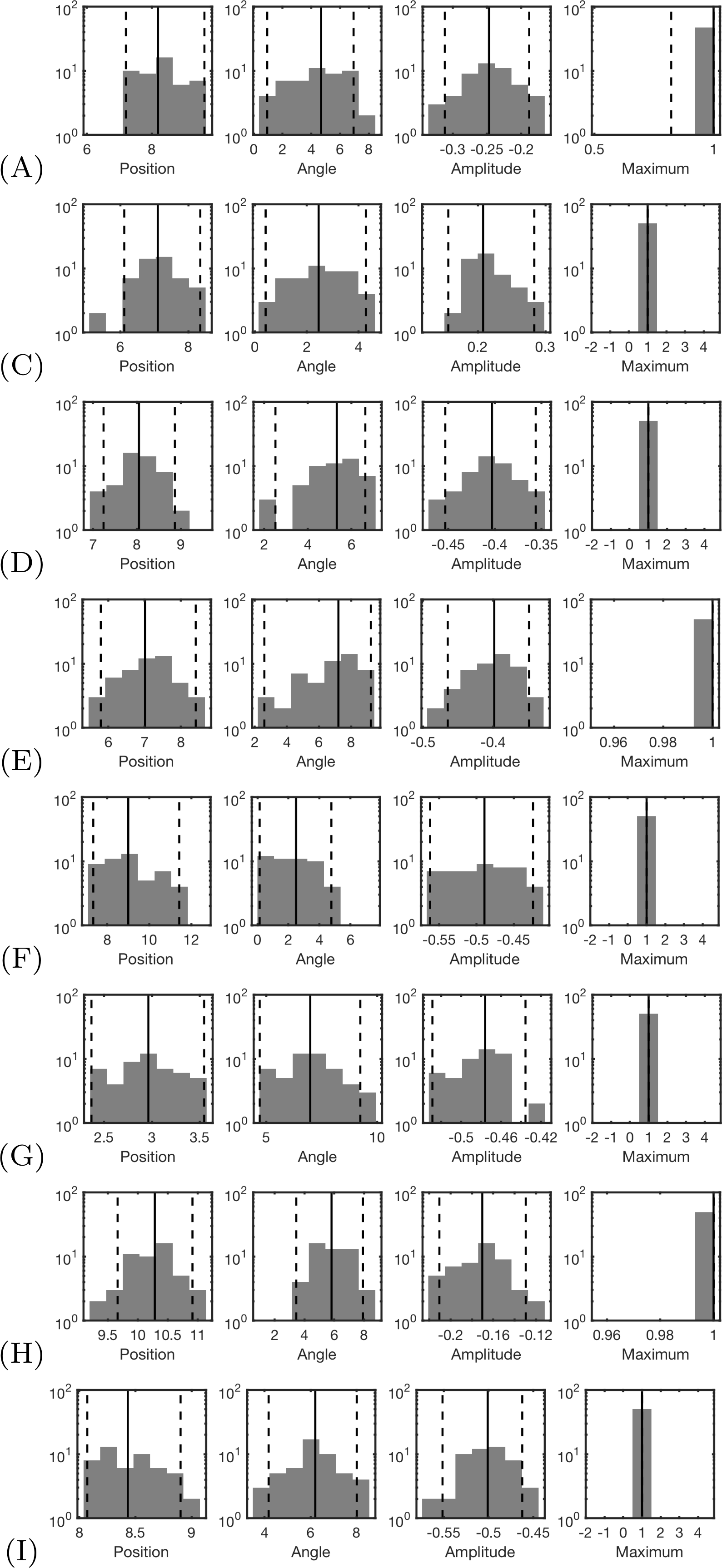}
    \caption{The results of the superficial source localization for the numerical experiments (A)--(I) conducted in the spherical domain (Table \ref{tab:A-H}). In contrast to the case of the deep source, the superficial one is detected accurately in each case where its amplitude differs from zero. The most accurate results were obtained, when the deep source was absent. The E/MEG yielded superior result compared to EEG. }
    \label{fig:histograms_superficial}
\end{figure}

\begin{figure}[h!]
    \centering
    \includegraphics[width=6cm]{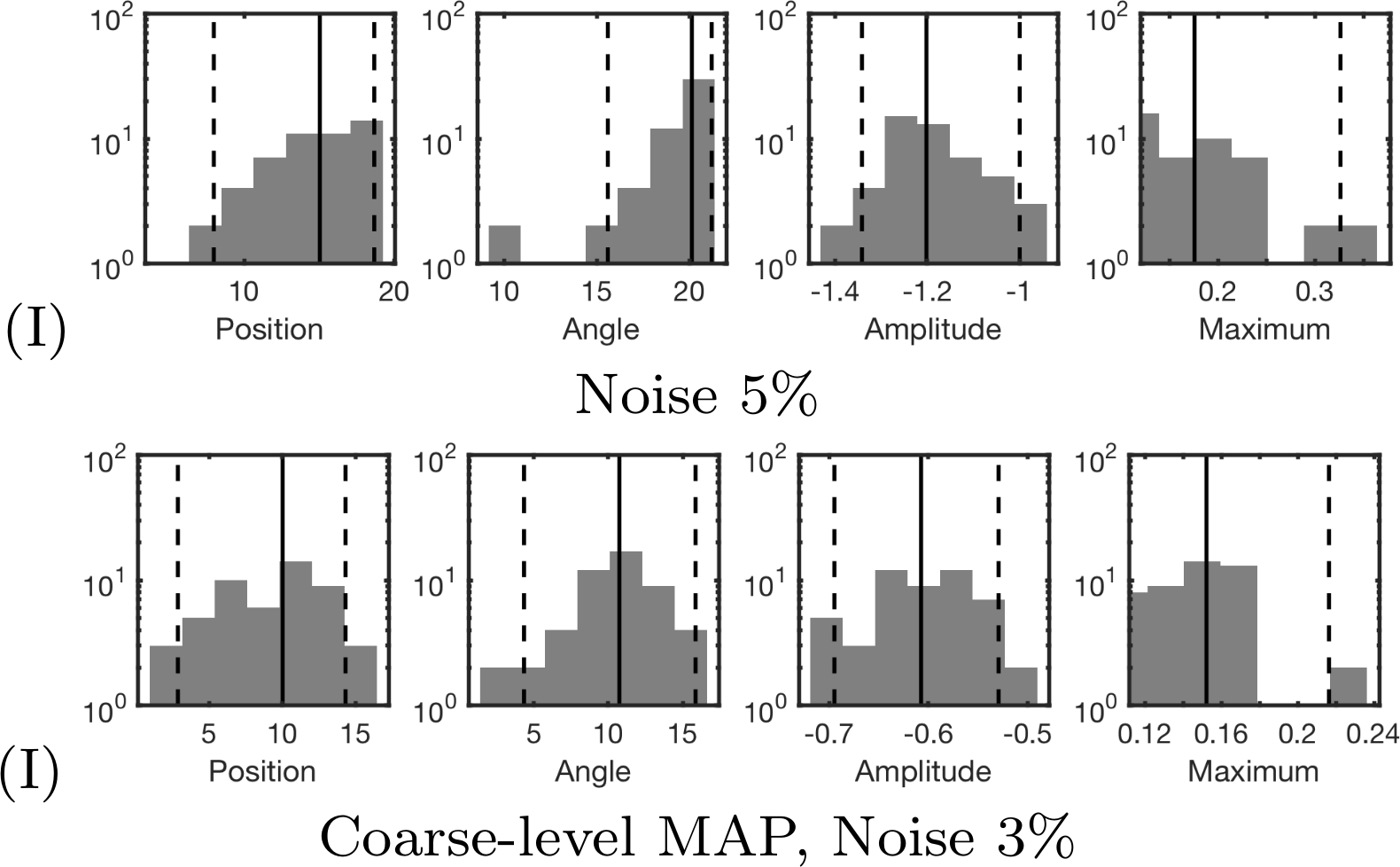}
    \caption{The additional deep source localization results obtained in the case (I) for the deep source. }
    \label{fig:histograms_deep_2}
\end{figure}

\begin{figure}[h!]
    \centering
   \includegraphics[width=6.5cm]{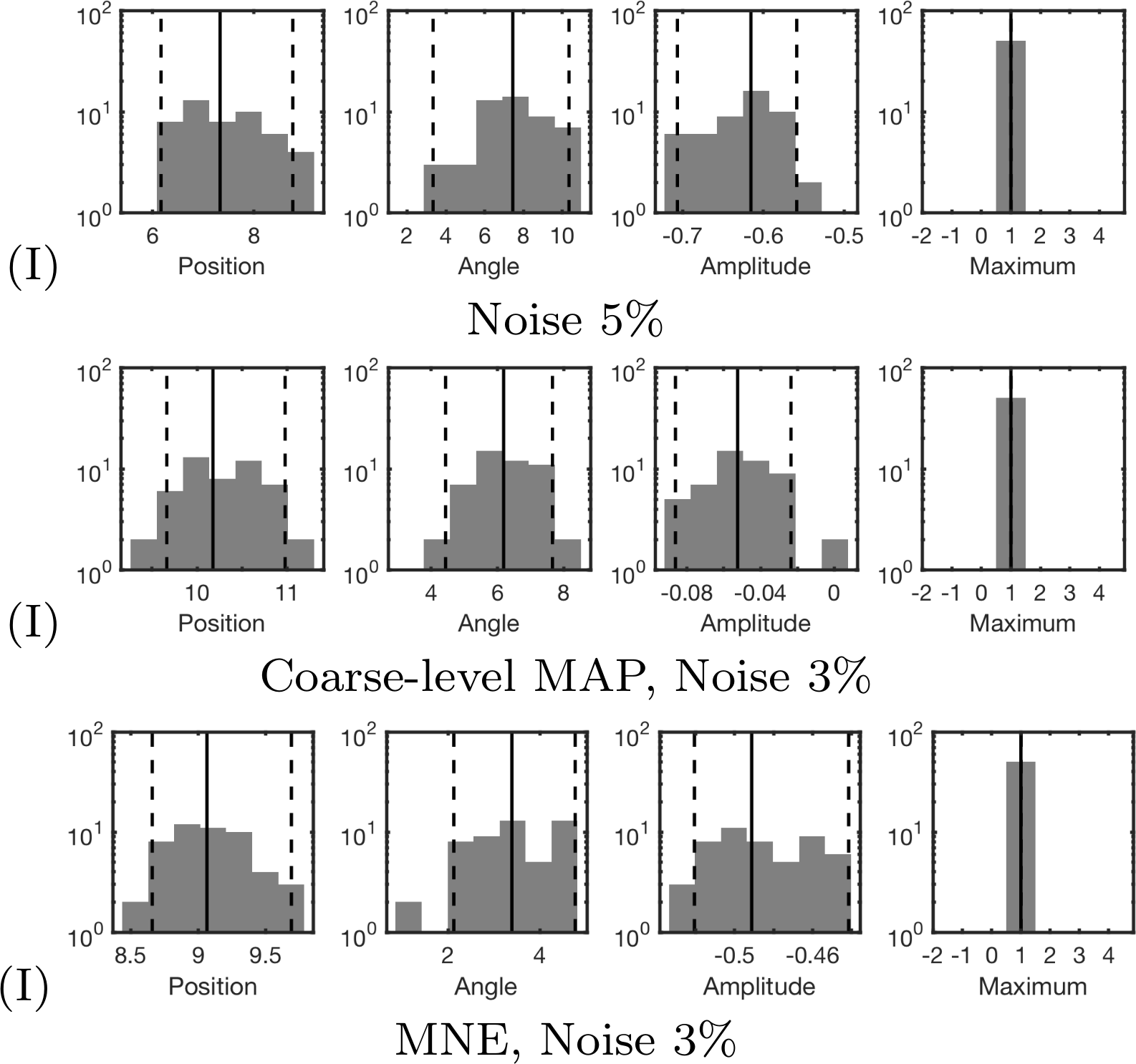}
    \caption{The additional superficial source localization results obtained in the case (I) for the superficial  source.}
    \label{fig:histograms_superficial_2}
\end{figure}

\begin{figure}[h]
     \centering
    \includegraphics[width=5cm]{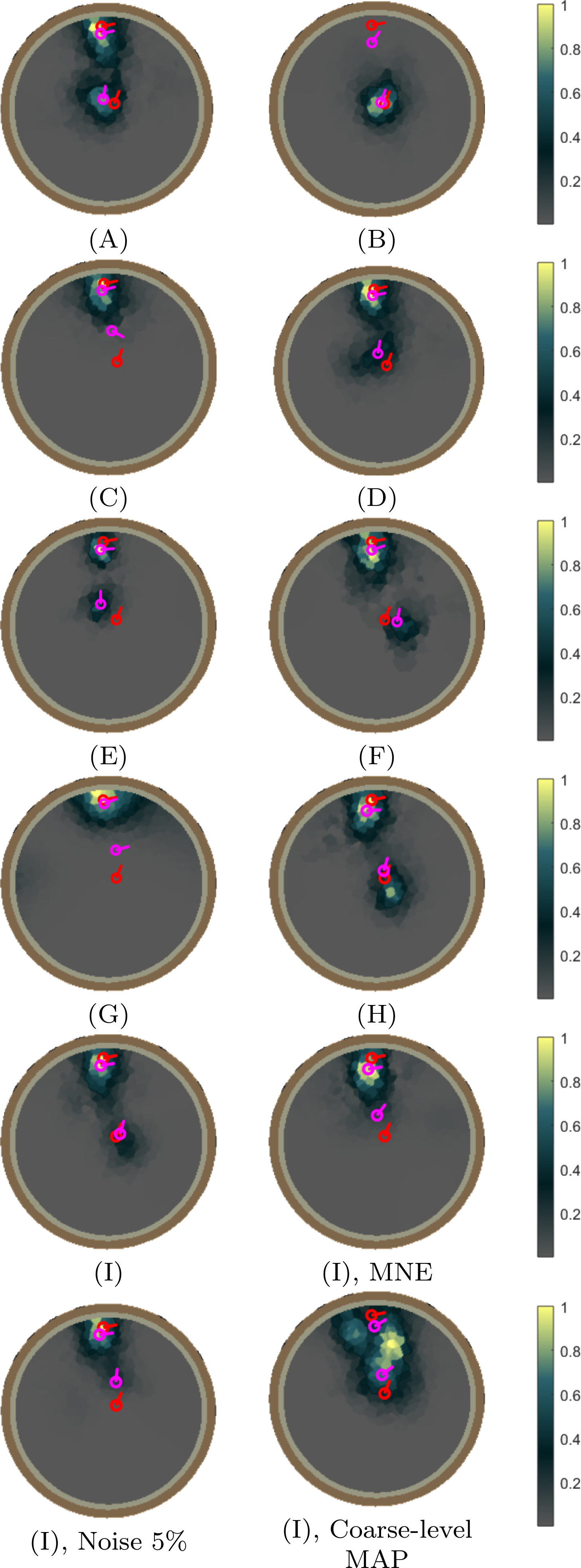}
    \caption{Examples of the reconstructions obtained in the  numerical experiments (A)--(I) in the spherical domain (Table \ref{tab:A-H}). In each image, the actual source and the center of mass of the reconstruction w.r.t.\ a ROI centered at the actual source position, are marked by the red and purple arrow, respectively.The case (I), was  studied using the following three alternative approaches in addition to the basic multiresolution scheme. MNE: only single IAS MAP iteration was performed on each reconstruction level, meaning that the estimate obtained coincided with MNE. Noise 5\%: the noise level was increased to 5\%. Coarse-level MAP: only the coarse resolution level was applied in the MAP estimation process with otherwise unchanged parameters. }
    \label{fig:spherical}
    \end{figure}

\begin{figure}[h!]
          \centering 
    \includegraphics[width= 6cm]{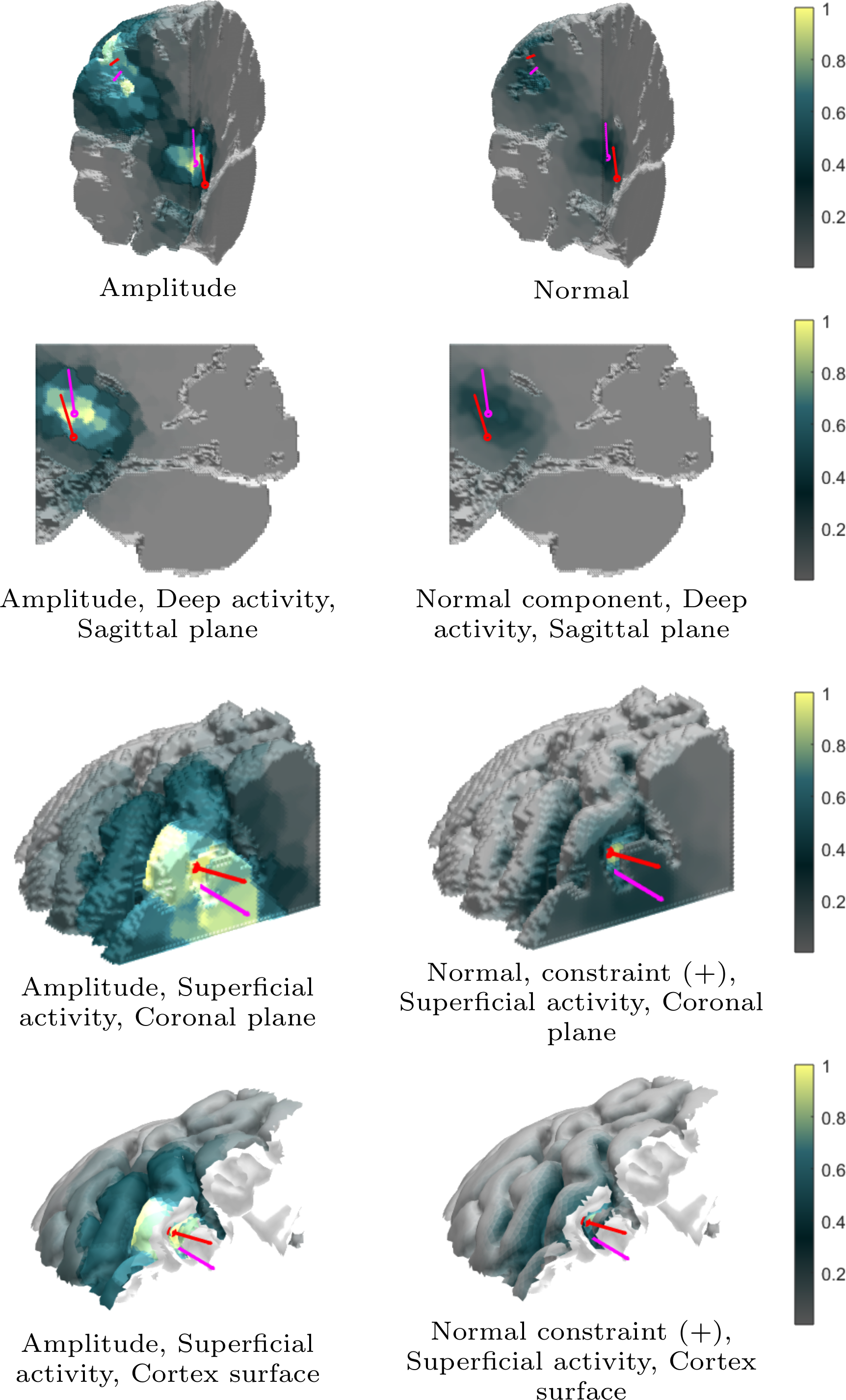} 
         \caption{The reconstruction (I) of the primary current density for the numerically modeled deep (thalamic P14/N14) and superficial (somatosensory P20/N20) activity obtained using the population head model (PHM).  On each row, the left column shows the amplitude and the right one the normal component in the direction of the surface normal. On 3rd and 4th row, the normal activity has been constrained into the outward direction. In each image, the actual source and the center of mass of the reconstruction w.r.t.\ a ROI centered at the actual source position, are marked by the red and purple arrow, respectively.  }
         \label{fig:realistic_volume}
     \label{fig:realistic_surface}
    \end{figure}

The results obtained in the numerical experiments have been included in Table \ref{tab:computing_time} and \ref{tab:my_label} and Figures \ref{fig:histograms_deep}--\ref{fig:realistic_surface}. In each case,  the deep and superficial component have been analyzed separately. Histograms for the cases (A)--(I) illustrate the accuracy of the reconstructed source with respect to the source position (mm), orientation (deg), amplitude, and the relative maximum of the current density within the ROI (Figures \ref{fig:histograms_deep}--\ref{fig:histograms_superficial_2}) with respect to the global maximum. The last one of these is utilized as a measure for the distinguishability of the source within the ROI. Examples of the reconstructions in the cases (A)--(I) are illustrated in Figure \ref{fig:spherical}, and the distributions obtained for (J), i.e., the realistic PHM,  are shown in Figure  \ref{fig:realistic_volume}. The additional cases evaluated for (I), are presented in Figure \ref{fig:histograms_deep_2} and \ref{fig:histograms_superficial_2}. 

The histograms in Figures \ref{fig:histograms_deep}--\ref{fig:histograms_superficial_2} illustrate the numerical accuracy of the RAMUS reconstruction  approach. Case (A) suggests that the activity in both superficial and deep areas can be reconstructed in EEG, when applying IG as hyperprior. In (A), the superficial source is found with the median positioning accuracy of 8 mm, angle difference of 4.5 degrees and  logarithmic (log10) relative amplitude error of -0.25, i.e., the amplitude of the reconstructed source is 56\% compared to that of the actual one. For the deep source these errors are 15 mm, 12 deg, -0.65 (22\% amplitude), respectively. Furthermore, as shown by the relative maximum, the superficial source always maximizes the (global) reconstruction, and the relative maximum within the deep ROI is around 50\% of the global one in median.

Based on (B) and (C), it is obvious that the reconstruction accuracy is better, if only one of the two sources is active. Furthermore, increasing the intensity of the superficial source decreases the reconstruction accuracy for deep one which is shown by the case of (D) for which the median position, orientation amplitude, and relative maximum for the deep source are 18 mm, 17 degrees, -0.85 and 0.25, respectively. That is, the accuracy is lower than in (A). In (E), a sparsity factor of 5 was used instead of 8, meaning that the resolution difference between the subsequent levels was less steep, resulting in a weaker distinguishability of the deep source. The same observation was made in the case (F) in which 20 randomized decompositions instead of 100 were used. The deep activity was absent in (G), where we used the G hyperprior instead of IG, confirming the necessity of IG as the hyperprior. In (H) and (I), the use of the E/MEG lead field was observed to improve the deep localization accuracy around 7 mm and orientation accuracy about 8 degrees with respect to the corresponding cases (A) and (D) of EEG data, while the superficial localization  accuracy was practically unchanged for (H) and deviated less than 2 mm and 1 deg for (I). The results for E/MEG were visually more focal than the ones obtained with EEG (Figure  \ref{fig:spherical}). 

In the three additional tests performed with the parameter setting (I), the increased 5\%  noise level led to 5 mm and 3 deg lower positioning and orientation accuracy for the deep source, and a smaller 1 mm and 1 deg deviation for the superficial one.  
In the case of the coarse-level MAP iteration with 3\% noise, 2 mm and 1 deg position and orientation improvement was observed for the deep source. For the superficial one, there was a 2 mm deviation in the position, while the orientation accuracy remained unchanged. The coarse-level estimate was visually less focal compared to the ones obtained with multiple resolution levels. MNE detected only the superficial source for which 1 mm position deflection and 3 deg orientation improvement were obtained compared to the basic case (I).


    \mbox{}

In (J), simultaneous localization of the simulated radial thalamic and tangential somatosensory component was found to be feasible with the realistic PHM model (Figure \ref{fig:realistic_volume}). 
Similar to the spherical case, the deep activity had a lower amplitude than the superficial one. In the somatosensory area, with the physiological normal constraint, i.e., the assumption that the primary current is oriented along the inward surface normal, the activity was localized in the posterior wall of the central sulcus similar to the synthetic P20/N20 component used in generating the data. 

The results concerning the computing times have been included in Table  \ref{tab:computing_time}. Those show that a superior  performance was obtained with GPU processing which provided the randomized set of decompositions and a  reconstruction in 1/5 to 1/4 of the time required by the laptop CPU.
\mbox{} 

\begin{table}[!h]
    \centering
        \caption{The computing time (in seconds) for 100 random  three-level multiresolution decompositions and of a  corresponding RAMUS (randomized multiresolution scan) estimate obtained with a NVIDIA Quadro P6000 workstation GPU and Intel i7 5650U laptop CPU.}
    \label{tab:computing_time}
    \begin{footnotesize}
    \begin{tabular}{l l l l l l l}
 Processor & Dec. &  EEG & E/MEG  \\ 
 \hline 
Quadro P6000 & 36 & 14 & 28 \\ 
i7 5650U  & 176 & 55 & 116 \\
    \end{tabular}
    \end{footnotesize}
\end{table}
\begin{table*}[h]
    \begin{center}
        \caption{The percentage of the reconstructions  fulfilling  the source detection criterion (the local maximum in the ROI \textgreater0.1 of the global maximum). This threshold criterion was chosen as it represents roughly the limit of a visually detectable source. In case (B), where only the deep source is present there are 4\% false positive detections for the superficial one. This is due to the relatively large deviation of the deep source due to which it is moved partially in the superficial ROI in the corresponding estimates. }
    \label{tab:my_label}
    \begin{tabular}{llllrrrrrrrrr}
     Type &  {\bf n}$\mbox{}^{a}$  & $\ell\mbox{}^{b}$ & ROI & (A) & (B) & (C) & (D) &  (E) & (F) & (G) & (H) & (I) \\ 
      \hline 
       MAP & 3 & 1--3 & Deep  & 100 & 100 &  0 & 98  &  100 & 86 & 0 & 100 & 92 \\
       & & & Sup. & 100 & 4 & 100 & 100 &  100 & 100 & 100 & 100 &  100 \\
        & 5 &  & Deep  &  &  &   &   &  &  &  & 100 & 90 \\
       & & & Sup. &  &  &  &  &   & &  & 100 &  100 \\
         MNE    & 3 &  1--3 & Deep  &  &  &  & &   &  &  &  & 0   \\
       & & & Sup. &  &  &  &  &   & &  & & 100    \\
                MAP    & 3 &  1 & Deep  &  &  & &   &   &  &  &  & 96             \\
       & & & Sup. &  &  &  &  &   & &  & & 100    \\
    \end{tabular} 
    \end{center}
    $\mbox{}^{a}$ Noise level (\%). \\
    $\mbox{}^{b}$ The resolution levels ($\ell$) in the reconstruction process.    
\end{table*}

\label{results}

\section{Discussion}
\label{discussion}

The present numerical results suggest that via the proposed randomized multiresolution scanning (RAMUS) technique one can obtain a robust and accurate MAP estimate for the primary current density in both superficial and deep parts of the brain. RAMUS was observed to enhance the visibility of deep components and also marginalizing the effect of the discretization without a remarkable  computation cost. The noise-robustness of RAMUS was shown for 3\% and 5\% noise levels. As expected, the effect of the noise was observed to be the most obvious with respect to the deep source. 

Utilizing a multiresolution approach was found to be crucial {\em per se} for the reconstruction quality, since maximal achievable accuracy for the deep components is significantly lower than for the superficial one. Detecting the deep source  necessitated the presence of a coarse resolution level in the MAP estimation process, i.e., a sparsity factor $s$ larger than one. The superior results were obtained with $s = 8$. Decreasing the value of $s$, i.e., increasing the source count, quickly diminished the detectability of the deep component which can be observed based on the results obtained for $s=5$.  The distinguishability of the deep source in the final estimate was  determined by the number of the source positions at the coarsest level which, in this study, was observed to be around 100--400 roughly matching the sparsity factors between $s = 8$ and $s=5$. Investigating this interval was motivated by the fact that the maximal number of the detectable sources in the numerical system is determined by the number of the data entries (Section \ref{sec:coarse_to_fine}) which is 102 for EEG and 204 for E/MEG, i.e., roughly of the same magnitude. In practice, the optimal size of the coarse system should also take the physiological modeling aspects into account and might be, therefore, also considerably larger than the present choice. 
For example, if the neural activity is limited to {\em a priori} known ROIs a larger number might be well-motivated.  A comparison between the single (coarse-level) and multiple resolution results showed that the refinement of the resolution during the reconstruction process improves the focality of the reconstruction and its accuracy in the superficial areas. Nevertheless, the coarse-level reconstruction was marginally superior in the deep part, emphasizing that here the finer resolution levels slightly affected the coarser level outcome, which is here presumably  optimal  for the weakly distinguishable deep activity. Thus, it is important to adjust the decomposition parameters appropriately. The marginalization of the discretization errors via random scanning was perceived to be vital in order to optimize the robustness of the reconstruction which was observed to grow along the number of the multiresolution decompositions utilized.

When coupled with the iterative alternating sequential (IAS) algorithm, RAMUS constitutes essentially a repetitive MAP optimization process for HBM  \citep{ohagan2004,calvetti2009}. Marginalizing the result over a given number of random multiresolution decompositions can be associated with computing an equal number of MAP estimates. Since the computational cost of the IAS algorithm is largely determined by the product between the lead field matrix and  a candidate solution which is  parallelized effectively in both CPU and GPU processors. Here, the latter option was found to achieve the fastest performance with the total computation time for a single reconstruction being 14 seconds which would be feasible in processing a larger dataset. Overall, the computational effort of evaluating the MAP via the RAMUS technique can be regarded as moderate compared to a full MCMC sampling based conditional mean (CM) estimate for the poerior density which has been evaluated in \citep{calvetti2009} within a ROI. Namely, achieving a full  convergence of MCMC would require  thousands of iterations \citep{liu2001} and the effort of one  iteration step is comparable to a single step of IAS. Thus, MCMC would be a slower option. Even though an optimization method, RAMUS can be also interpreted as a surrogate  for CM, as it, on one hand,  increases the robustness of the source reconstruction via sampling, but, on the other hand, does not provide  as extensive information about the posterior density itself as an actual Bayesian sampler does. 

 The results obtained suggest that the IG hyperprior \citep{ohagan2004} is necessary in conjunction with IAS, when it is coupled with RAMUS, as the deep activity was not detected with G. Since here the cases of the G hyperprior and single-step MAP can be associated with the 1-norm regularized MCE and MNE \citep{uuhaso99,hamalainen1993}  (Section \ref{sec:hbm}), respectively, it also seems that RAMUS would not enable correcting a depth bias related to either of these estimates. Previously, in  \citep{calvetti2009}, IG was found to perform well for the deep part, when a region of interest was used. Based on the present results, RAMUS provides the means to utilize the advantage of the IG within the whole brain and with a high imaging resolution, while maintaining the computational cost on an appropriate level. 
 
 We emphasize that the present conditionally Gaussian prior, in its current formulation, is depth, resolution and decomposition invariant.  That is, additional physiological or operator based weighting or prior conditioning  \citep{homa2013bayesian,calvetti2015hierarchical,calvetti2018brain} is not necessary in order to balance the depth performance of the MAP estimate. Our interpretation for this is that RAMUS can correct the depth localization inaccuracies that are otherwise found with MAP estimates, as it, via the multiresolution approach, decomposes the source space into a  set of a visible and invisible fluctuations, explores both sets, and also helps to marginalize the random numerical discretization and optimization errors out of the final estimate. Central here are the visibility of the deep activity at low resolution levels \citep{pascual1999review}, the concept of the sensitivity decomposition \citep{liu1995sensitivity} and forming such through  projections and multiresolution decompositions which have been investigated in the context of other inverse problems, e.g., in \citep{piana1997projected,pursiainen2008_3}. 
 
 In addition to the investigated properties, the choice for the scale parameter was also  observed to be important in order to guarantee proper function of RAMUS. In each MAP estimation process, the present value 1E-10 was found to work well in detecting activity for both the spherical and realistic head model. The workable  range for the scale parameter was observed to be from 1E-10 to 1E-08 similar to the previous findings \citep{calvetti2009}. Outside this interval the deep activity was not found appropriately or the orientation accuracy of the estimates was lost. In the latter case, the estimate was locked into the direction of a Cartesian coordinate axis, meaning that, due to overly strong focality condition, only single component in the estimate differed from zero in the end of the iteration. Locking was also observed for MAP optimization sequences considerably longer than 10 iteration steps. With a sufficiently large scale parameter there was no locking, but the reconstructions obtained  were also less focal.
 
 The results of this article concern only the present numerical framework in which a deep and superficial source were detected simultaneously. Future work will include testing and analyzing the performance of the RAMUS approach with real experimental SEP/F data with the goal to distinguish cortical and sub-cortical activity, e.g., the P14/N14 (deep)  and P20/N20 (superficial) components occuring in the stimulation of the median nerve. A comparison with other inverse methods capable of  deep localization, such as LORETA and Beamforming \citep{pascual2002functional,pascual1999low,jonmohamadi2014comparison},  will also be important. Further method development topics will include a deeper  investigation on the inverse effects of the hyperprior and decomposition parameters as well as finding alternative strategies to update the initial guess for the IAS MAP estimation technique. In the latter case, for instance, an averaged initial guess obtained with respect to several multiresolution decompositions might provide a potential alternative for the current approach which relies on a single decomposition. To make the random scanning computationally more efficient a solver based on parallel scanning processes might be developed. We also  emphasize that RAMUS with its current formulation, the proposed algorithm can be applied to reduce discretization errors not only with the present IAS MAP method but potentially for a variety  of source reconstruction techniques.
 

\section*{Acknowledgments}

AR and SP were supported by the Academy of Finland Centre of Excellence in Inverse Modelling and Imaging 2018--2025 and the Vilho, Yrjö and Kalle Väisälä Foundation of Finnish Academy of Science and Letters. AK was supported by the Academy of Finland Postdoctoral Researcher grant number 316542.

\appendix

\section{IAS MAP estimation}
\label{sec:ias}

The IAS algorithm finds a MAP estimate for the posterior $p ( {\bf x}, {\bf \bm \theta} \mid {\bf y})$ as follows:
 \begin{enumerate}
\item Set $k=0$ and ${\bm \theta}^{(0)} = (\theta_0, \theta_0, \ldots, \theta_0)$.
\item Set ${\bf L}^{(k)} = {\bf L} {\bf D}^{1/2}_{{\bm \theta}^{(k)}}$ with \begin{equation} {\bf D}^{1/2}_{{\bm \theta}^{(k)}}= \hbox{diag} (\sqrt{ |{\bm \theta}^{(k)}_1|}, \sqrt{|{\bm \theta}^{(k)}_2|}, \ldots, \sqrt{|{\bm \theta}^{(k)}_n|}). \end{equation}
\item Evaluate \begin{equation}
{\bf x}^{(k+1)} =  {\bf D}^{1/2}_{{\bm \theta}^{(k)}} {{\bf L}^{(k)}}^T ( {\bf L}^{(k)} {{\bf L}^{(k)}}^T + \sigma^2  {\bf I})^{-1} {\bf y},  \label{ias_matrix} \end{equation}
 where $\sigma$ denotes the standard deviation of the likelihood.
\item Update the hyperparameter based on the hypermodel. 
\begin{itemize} 
\item If the hypermodel is G, set \begin{equation} \theta_i = \frac{ 1 } { 2 }  \theta_0 \Big( \eta + \sqrt{\eta^2 + 2 {x_i^{(k)}}^2/\theta_0} \Big) \end{equation} with $\eta = \beta - 3/2$, $i = 1, 2, \ldots, n$.
\item Else, if the hypermodel is IG, set \begin{equation}\theta_i^{(k+1)} =  (\theta_0 +  \frac{{x_i^{(k)}}^2}{2})/\kappa \end{equation}  with $\kappa = \beta + 3/2$, $i = 1, 2, \ldots, n$.
\end{itemize}
\item Set $k = k +1$ and go back to 2., if $k$ is less than the total number of iterations defined by the user. \end{enumerate}

\end{document}